\documentclass[11pt]{article}

\usepackage{amssymb,latexsym}
\usepackage{epsfig}
\usepackage{eufrak}
\usepackage{amsmath}
\usepackage{mathrsfs}
\usepackage{color}

\newtheorem{theorem}{Theorem}
\newtheorem{lemma}{Lemma}

\newcommand{\be}{\begin{equation}}
\newcommand{\ee}{\end{equation}}
\newcommand{\bee}{\begin{eqnarray*}}
\newcommand{\eee}{\end{eqnarray*}}
\newcommand{\bel}{\begin{eqnarray}}
\newcommand{\eel}{\end{eqnarray}}
\newcommand{\bec}{\begin{cases}}
\newcommand{\eec}{\end{cases}}
\newcommand{\bem}{\begin{bmatrix}}
\newcommand{\eem}{\end{bmatrix}}

\newcommand{\la}{\label}
\newcommand{\li}{\left}
\newcommand{\ri}{\right}

\newcommand{\ovl}{\overline}

\newcommand{\lc}{\lceil}
\newcommand{\rc}{\rceil}
\newcommand{\lf}{\lfloor}
\newcommand{\rf}{\rfloor}

\newcommand{\lm}{\lambda}

\newcommand{\vDe}{\varDelta}

\newcommand{\al}{\alpha}
\newcommand{\ba}{\beta}

\newcommand{\Om}{\Omega}

\newcommand{\f}{\frac}
\newcommand{\sq}{\sqrt}
\newcommand{\cd}{\cdots}

\newcommand{\qu}{\quad}
\newcommand{\qqu}{\qquad}
\newcommand{\fa}{\forall}

\newcommand{\mscr}{\mathscr}

\newcommand{\bb}{\mathbb}

\newcommand{\mrm}{\mathrm}

\newcommand{\tx}{\text}

\newcommand{\iy}{\infty}

\newcommand{\pa}{\partial}

\newcommand{\bed}{\begin{description}}
\newcommand{\eed}{\end{description}}
\newcommand{\bei}{\begin{itemize}}
\newcommand{\eei}{\end{itemize}}
\newcommand{\ben}{\begin{enumerate}}
\newcommand{\een}{\end{enumerate}}
\newcommand{\bib}{\bibitem}
\newcommand{\beL}{\begin{lemma}}
\newcommand{\eeL}{\end{lemma}}
\newcommand{\beT}{\begin{theorem}}
\newcommand{\eeT}{\end{theorem}}
\newcommand{\sect}{\section}

\newcommand{\bpf}{\begin{pf}}
\newcommand{\epf}{\end{pf}}
\newcommand{\bsk}{\bigskip}
\newcommand{\bi}{\binom}

\newcommand{\pfbox}{\hfill\mbox{$\Box$}}

\newenvironment{pf}{\paragraph*{Proof{\rm.}}}{\pfbox\bigskip}

\begin{document}

\title{{\bf A Truncation Approach for Fast Computation of Distribution Functions}
\thanks{The author is currently with Department of Electrical Engineering,
Louisiana State University at Baton Rouge, LA 70803, USA, and Department of Electrical Engineering, Southern
University and A\&M College, Baton Rouge, LA 70813, USA; Email: chenxinjia@gmail.com}}

\author{Xinjia Chen}

\date{February, 2008}

\maketitle

\begin{abstract}

In this paper, we propose a general approach for improving the efficiency of computing distribution functions. The idea is to truncate the
domain of summation or integration.

\end{abstract}

\sect{General Theory}

In various fields of sciences and engineering, it is a frequent problem to compute distribution functions.  Specifically, it is desirable to
compute efficiently and precisely the probability
\[
P = \Pr \{  a_i < X_i < b_i, \; i = 1 , \cd, m   \}
\]
where $X_1, \cd, X_m$ are random variables defined in probability space $(\Om, \mscr{F}, \Pr)$ and $a_i, \; b_i, \; i = 1, \cd, m$ are real
numbers.  Since the probability can be expressed as an $m$-dimensional summation or integration over domain
\[
D = \{(x_1, \cd, x_m): a_i < x_i < b_i, \; i = 1 , \cd, m \},
\]
the computational complexity may depend on the size of $D$.  In many situations, the larger the size of domain $D$ is, the more computation is
required. Hence, it will be of computational advantage to reduce the domain $D$ as its subset
\[
D^\prime = \{(x_1, \cd, x_m): a_i^\prime < x_i < b_i^\prime, \; i = 1 , \cd, m \}
\]
such that
\[
P^\prime = \Pr \{ a_i^\prime < X_i < b_i^\prime, \; i = 1 , \cd, m \}
\]
is close to $P$ within a controllable difference. For this purpose, we have

\beT Let $u_i, \; v_i, \; \al_i, \; \ba_i$ be real numbers such that
\[
\Pr \{ X_i \leq u_i \} \leq \al_i, \qqu \Pr \{ X_i \geq v_i \} \leq
\beta_i
\]
for $i = 1, \cd, m$.  Let
\[
a_i^\prime = \max(a_i, u_i), \qqu b_i^\prime = \min(b_i, v_i), \qqu i = 1 , \cd, m.
\]
Then, \[ P^\prime \leq P \leq P^\prime + \sum_{i = 1}^m (\al_i + \ba_i).
\]
\eeT

\bpf Obviously, $P^\prime \leq P$ is true since $D^\prime$ is a subset of $D$. Thus, it suffices to show $P \leq P^\prime + \sum_{i = 1}^m
(\al_i + \ba_i)$.

 Note that \[ P^\prime  =  \Pr \{ a_i^\prime < X_i < b_i^\prime, \; i = 1 , \cd, m \} =  \Pr \{ \cap_{i = 1}^m (A_i \cap B_i)
\cap C \}
\] where
\[
A_i = \{ X_i > u_i \}, \qqu B_i = \{ X_i < v_i \}, \qqu  i = 1 ,
\cd, m
\]
and $C = \{  a_i < X_i < b_i,  \; i = 1 , \cd, m \}$.   By Bonferroni's inequality,
\[
\Pr \{ \cap_{i = 1}^m (A_i \cap B_i) \cap C \} \geq  \Pr \{C\} - 2m + \sum_{i = 1}^m ( \Pr \{A_i\} + \Pr \{B_i \} ).
\]
By the definitions of $P$ and $P^\prime$,
\[
P^\prime \geq P - 2m + \sum_{i = 1}^m ( \Pr \{A_i\} + \Pr \{B_i \} ).
\]
Hence, \bee
 P & \leq & P^\prime + 2m -  \sum_{i = 1}^m ( \Pr \{A_i\} + \Pr \{B_i \} )\\
 & = & P^\prime +  \sum_{i = 1}^m \li [ ( 1 - \Pr \{A_i\}) + (1 -  \Pr \{B_i \} ) \ri ]\\
 & = & P^\prime +  \sum_{i = 1}^m \li [ \Pr \{X_i \leq u_i\} + \Pr \{X_i \geq v_i\} \ri ]\\
 & \leq & P^\prime + \sum_{i = 1}^m (\al_i + \ba_i).
\eee

This completes the proof of the theorem.

\epf

To ensure that $P^\prime \leq P \leq P^\prime + \eta$ for a prescribed $\eta > 0$, it suffices to choose
\[
\al_i = \ba_i = \f{\eta}{2 m}.
\]

As can be seen from Theorem 1, a critical step is to determine $u$ and $v$ for a random variable $X$ such that
\[
\Pr \{ X \leq u \} \leq \al, \qqu \Pr \{ X \geq v \} \leq \ba
\]
for prescribed $\al, \; \ba \in (0,1)$.  For this purpose, we have the following theorem.

\beT

Let $X$ be a random variable with mean $\mu = \bb{E} [X]$.  Let
\[
\mscr{C} (z) = \inf_{t \in \bb{R}} \bb{E} [ e^{t (X - z)} ] .
\]
Then the following statements hold true:

(I) For any $z > \mu$,
\[
\Pr \li \{  X  \geq z \ri \} \leq \mscr{C} (z).
\]

(II) For any $z < \mu$,
\[
\Pr \li \{ X  \leq z \ri \} \leq \mscr{C} (z).
\]

(III) Both $\mscr{C} (\mu + \vDe)$ and $\mscr{C} (\mu - \vDe)$ are monotonically decreasing with respect to $\vDe > 0$.

(IV) For any $\al \in (0,1)$, there exists a unique $\vDe > 0$ such that $\mscr{C} (\mu - \vDe) = \al$.

(V) For any $\ba \in (0,1)$, there exists a unique $\vDe > 0$ such that $\mscr{C} (\mu + \vDe) = \ba$.

 \eeT

\bpf By Jensen's inequality
\[
\bb{E} [ e^{t (X - z)} ] \geq  e^{t \bb{E}[X - z]}.
\]
Hence, if $z < \mu$, we have $\bb{E} [ e^{t (X - z)} ] \geq  e^{t \bb{E}[X - z]} \geq 1$ for $t  \geq 0$.  Similarly, if $z > \mu$, we have
$\bb{E} [ e^{t (X - z)} ] \geq  e^{t \bb{E}[X - z]} \geq 1$ for $t  \leq 0$.  Combing these observations and the fact that
\[
0 \leq \mscr{C} (z) \leq \bb{E} [ e^{0 \times (X - z)} ] = 1,
\]
we have
\[
\mscr{C} (z) =  \bec \inf_{t < 0}   \bb{E} [
e^{t (X - z)} ]  & \tx{for} \; z < \mu,\\
\inf_{t > 0}  \bb{E} [ e^{t (X - z)} ]  & \tx{for} \; z > \mu. \eec
\]
By the Chernoff bounds \cite{Chernoff},
\[
\Pr \li \{  X  \geq z \ri \} \leq \inf_{t < 0}   \bb{E} [ e^{t (X - z)} ]   = \mscr{C} (z) \] for $z < \mu$; and \[ \Pr \li \{  X  \leq z \ri \}
\leq \inf_{t > 0}   \bb{E} [ e^{t (X - z)} ]   = \mscr{C} (z) \] for $z > \mu$.   This completes the proof of statements (I) and (II).

To show that $\mscr{C} (\mu + \vDe)$ is monotonically decreasing with respect to $\vDe > 0$,  let $t_\vDe$ be the number such that
\[
\inf_{t \in \bb{R}} \bb{E} [ e^{t (X- \mu - \vDe)}] = \bb{E} [ e^{t_\vDe (X- \mu - \vDe)}].
\]
Then, $t_\vDe$ is positive and
\[
\f{\pa    \bb{E} [ e^{t_\vDe (X- \mu + \vDe)}] } {\pa t_\vDe } = 0.
\]
It follows that \bee \f{d   \mscr{C} (\mu + \vDe)} {d \vDe  } & = & \f{d \inf_{t \in \bb{R}}
\bb{E} [ e^{t (X- \mu - \vDe)}] } {d \vDe  }\\
& = & \f{d  \; \bb{E} [ e^{t_\vDe (X- \mu - \vDe)}] } {d
\vDe  } \\
& = & \f{\pa    \bb{E} [ e^{t_\vDe (X- \mu - \vDe)}] } {\pa \vDe } + \f{\pa    \bb{E} [ e^{t_\vDe (X- \mu - \vDe)}] } {\pa t_\vDe } \f{\pa
t_\vDe}
{\pa \vDe }\\
& = & \f{\pa    \bb{E} [ e^{t_\vDe (X- \mu - \vDe)}] } {\pa \vDe } +
0 \times \f{\pa t_\vDe} {\pa \vDe }\\
& = & \f{\pa    \bb{E} [ e^{t_\vDe (X- \mu - \vDe)}] } {\pa \vDe }  < 0. \eee

Similarly, to show that $\mscr{C} (\mu - \vDe)$ is monotonically decreasing with respect to $\vDe > 0$, let $t_\vDe$ be the number such that
\[
\inf_{t \in \bb{R}} \bb{E} [ e^{t (X- \mu + \vDe)}] = \bb{E} [ e^{t_\vDe (X- \mu + \vDe)}].
\]
Then, $t_\vDe$ is negative and
\[
\f{\pa    \bb{E} [ e^{t_\vDe (X- \mu + \vDe)}] } {\pa t_\vDe } = 0.
\]
Consequently, \bee \f{d   \mscr{C} (\mu - \vDe)} {d \vDe  } & = & \f{d \inf_{t \in \bb{R}}
\bb{E} [ e^{t (X- \mu + \vDe)}] } {d \vDe  }\\
& = & \f{d   \; \bb{E} [ e^{t_\vDe (X- \mu + \vDe)}] } {d
\vDe  } \\
& = & \f{\pa    \bb{E} [ e^{t_\vDe (X- \mu + \vDe)}] } {\pa \vDe } + \f{\pa    \bb{E} [ e^{t_\vDe (X- \mu + \vDe)}] } {\pa t_\vDe } \f{\pa
t_\vDe}
{\pa \vDe }\\
& = & \f{\pa    \bb{E} [ e^{t_\vDe (X- \mu + \vDe)}] } {\pa \vDe } +
0 \times \f{\pa t_\vDe} {\pa \vDe }\\
& = & \f{\pa    \bb{E} [ e^{t_\vDe (X- \mu + \vDe)}] } {\pa \vDe }  < 0. \eee  This concludes the proof of statements (III).

To show statement (IV), note that {\small \be \la{c3}
 \liminf_{\vDe \to 0} \mscr{C} (\mu - \vDe) = \liminf_{\vDe \to 0} \; \inf_{t < 0}  \bb{E} [ e^{t (X
- \mu + \vDe)} ] \geq \liminf_{\vDe \to 0} \; \inf_{t < 0}  e^{ \bb{E} [t(X - \mu + \vDe)]}  =  \lim_{\vDe \to 0} \; \inf_{t < 0} e^{ \vDe t} =
1 \ee}
 and that \be \la{c4}
 \lim_{\vDe \to \iy} \mscr{C} (\mu -
\vDe) = 0 \ee as a result of
\[
\liminf_{\vDe \to \iy} \mscr{C} (\mu - \vDe) \geq 0 \]
 and
\[
\limsup_{\vDe \to \iy} \mscr{C} (\mu - \vDe) = \limsup_{\vDe \to \iy} \inf_{t < 0}  \bb{E} [ e^{t (X - \mu + \vDe)} ] \leq \limsup_{\vDe \to
\iy} \bb{E} [ e^{ - (X - \mu + \vDe)} ]  = \bb{E} [ e^{ - (X - \mu )} ] \lim_{\vDe \to \iy} e^{ - \vDe}  = 0.
\]
Hence, (IV) follows from (\ref{c3}), (\ref{c4}) and the fact that $\mscr{C} (\mu - \vDe)$ is monotonically decreasing with respect to $\vDe >
0$.

To show statement (V), note that {\small \be \la{c1}
 \liminf_{\vDe \to 0} \mscr{C} (\mu + \vDe) = \liminf_{\vDe \to 0} \; \inf_{t > 0}  \bb{E} [ e^{t (X
- \mu - \vDe)} ] \geq \liminf_{\vDe \to 0} \; \inf_{t > 0}  e^{ \bb{E} [t(X - \mu - \vDe)]}  =  \lim_{\vDe \to 0} \; \inf_{t > 0} e^{ - \vDe t}
= 1 \ee}
 and that \be \la{c2}
 \lim_{\vDe \to \iy} \mscr{C} (\mu +
\vDe) = 0 \ee as a result of
\[
\liminf_{\vDe \to \iy} \mscr{C} (\mu + \vDe) \geq 0 \]
 and
\[
\limsup_{\vDe \to \iy} \mscr{C} (\mu + \vDe) = \limsup_{\vDe \to \iy} \inf_{t > 0}  \bb{E} [ e^{t (X - \mu - \vDe)} ] \leq \limsup_{\vDe \to
\iy} \bb{E} [ e^{ (X - \mu - \vDe)} ]  = \bb{E} [ e^{ (X - \mu )} ] \lim_{\vDe \to \iy} e^{ - \vDe}  = 0.
\]
Hence, (V) follows from (\ref{c1}), (\ref{c2}) and the fact that $\mscr{C} (\mu + \vDe)$ is monotonically decreasing with respect to $\vDe > 0$.

\epf

\bsk

As can be seen from Theorem 2, since $\mscr{C} (\mu - \vDe)$ is monotonically decreasing with respect to $\vDe > 0$, we can determine $\vDe > 0$
such that $\mscr{C} (\mu - \vDe) = \al$ by a bisection search. Then, setting $u = \mu - \vDe$ yields $\Pr \{X \leq u \} \leq \al$ as desired.
Similarly, we can determine $\vDe > 0$ such that $\mscr{C} (\mu + \vDe) = \ba$ by a bisection search and set $v = \mu + \vDe$ to ensure $\Pr \{X
\geq v \} \leq \ba$.

\section{Applications}

The approach of reducing the domain $D$ to its subset $D^\prime$ is referred to as {\it truncation technique} in this paper.  By the Chebyshev's
inequality, it can be visualized that if the variances of $X_i$ are small, then the size of the truncated domain $D^\prime$ can be much smaller
than that of domain $D$, even though $\eta$ is extremely small.

For the truncation technique to be of practical use, it is desirable that functions $\mscr{C}(z)$ associated $X_i$ have closed form. This is
indeed the case for many important distributions.  For example, when $X$ is the average of i.i.d Bernoulli random variables $Y_1, \cd, Y_n$ such
that $\Pr \{ Y_i = 1 \} = p$ for $1 \leq i \leq n$, the Hoeffding's inequality \cite{{Hoeffding}} asserts that
\[
\Pr \{ X \geq z  \} \leq \mscr{C}(z), \qqu \fa z > p
\]
\[
\Pr \{ X \leq z  \} \leq \mscr{C}(z), \qqu \fa z < p
\]
where
\[
\mscr{C}(z) = \li [ \li ( \f{p} {z}  \ri )^{ z } \li ( \f{1- p} {1- z}  \ri )^{ 1- z } \ri ]^n.
\]

For another example, when $X$ is the average of i.i.d Poisson random variables $Y_1, \cd, Y_n$ such that $\bb{E} \{ Y_i = 1 \} = \lm$ for $1
\leq i \leq n$, it can be shown by the Chernoff bounds \cite{Chernoff} that
\[
\Pr \{ X \geq z  \} \leq \mscr{C}(z), \qqu \fa z > \lm
\]
\[
\Pr \{ X \leq z  \} \leq \mscr{C}(z), \qqu \fa z < \lm
\]
where
\[
\mscr{C}(z) = \li [ e^{- \lm} \li ( \f{ \lm e } { z } \ri )^z \ri ]^n.
\]
Similar truncation techniques can be developed for hypergeometric distribution, negative binomial distribution, hypergeometric waiting-time
distribution, etc.

In the case that simple and tight bounds of $\mscr{C}(z)$ are
available, it is convenient to use the bounds in the truncation of
$D$. In this regard, we have established the following result.

\beT Let $K$ be a binomial random variable such that {\small $\Pr \{
K = i \} = \bi{n}{i} p^i (1 - p)^{n - i}, \; i = 0, 1, \cd , n$}
where $p \in (0, 1)$ and $n$ is a positive integer. Then, for
arbitrary real numbers $a, \; b$ and any $\eta \in (0,1)$, \[ \Pr \{
\mrm{T}^{-} \leq K \leq \mrm{T}^{+} \} \leq \Pr \{ a \leq K \leq b
\} < \Pr \{ \mrm{T}^{-} \leq K \leq \mrm{T}^{+} \} + \eta
\]
where {\small
\[
\mrm{T}^{-}  = \max \li \{ a, \; \li \lc n p + \f{ 1 - 2 p - \sq{ 1 + \f{18 n p (1-p)}{ \ln \f{2}{\eta}} } } {\f{2}{3n} + \f{3}{ \ln
\f{2}{\eta}}}  \ri \rc \ri \},  \qu \mrm{T}^{+} = \min \li \{ b, \; \li \lf n p + \f{ 1 - 2 p + \sq{ 1 + \f{18 n p (1-p)}{ \ln \f{2}{\eta}} } }
{\f{2}{3n} + \f{3}{\ln \f{2}{\eta}}} \ri \rf \ri \}
\]}
with $\lf . \rf$ and $\lc . \rc$ denoting the floor and ceiling functions respectively.
 \eeT

\bsk

We would like to remark that $\mrm{T}^{+} - \mrm{T}^{-}$ can be much
smaller than $b - a$ even though $\eta$ is chosen as an extremely
small positive number.

\bsk

To prove Theorem 3, we need some preliminary results.

\beL \la{lemcap} Define {\small $\mscr{M}(z, \mu) = \f{ (\mu - z)^2
} {2 \li ( \f{2 \mu}{3} + \f{z}{3} \ri ) \li ( \f{2 \mu}{3} +
\f{z}{3} - 1 \ri ) }$} for $0 < \mu < 1$ and $- 2 \mu < z < 3 - 2
\mu$. Then, for any fixed $\mu \in (0, 1)$,  $\mscr{M}(z, \mu)$ is
monotonically increasing from $- \iy$ to $0$ as $z$ increases from
$- 2 \mu$ to $\mu$, and is monotonically decreasing from $0$ to $-
\iy$ as $z$ increases from $\mu$ to $3 - 2 \mu$.

 \eeL

\bpf After a lengthy calculation, we obtained
\[
\f{\pa \mscr{M}(z,\mu)}{\pa z}  =  \f{(\mu - z) \; w(z, \mu)  } {
 \li [ \li ( \f{2\mu}{3} + \f{z}{3} \ri ) \li ( 1 - \f{2\mu}{3} - \f{z}{3} \ri ) \ri ]^2
 }
 \]
 where $w(z, \mu) = \mu(1 - \f{2\mu}{3} - \f{z}{3} ) + \f{z -
\mu}{6}$.  Noting that $w(- 2 \mu, \mu) = \f{\mu}{2} > 0, \; w( \mu,
\mu) = \mu (1 - \mu) > 0$ and that $w(z, \mu)$ is linear with
respect to $z$, we have that $w(z, \mu) > 0$ for any $\mu \in (0,
1)$ and $z \in (- 2 \mu, \mu)$.  It follows that $\f{\pa
\mscr{M}(z,\mu)}{\pa z}  > 0$ for any $\mu \in (0, 1)$ and $z \in (-
2 \mu, \mu)$.  In view of the positive sign of the partial
derivative and the fact that $\lim_{z \to - 2 \mu} \mscr{M}(z,\mu) =
- \iy, \; \mscr{M}(\mu,\mu) = 0$, we have that, for any fixed $\mu
\in (0, 1)$, $\mscr{M}(z, \mu)$ is monotonically increasing from $-
\iy$ to $0$ as $z$ increases from $- 2 \mu$ to $\mu$.

Similarly, observing that $w(3 - 2 \mu, \mu) = \f{1 - \mu}{2} > 0,
\; w( \mu, \mu) = \mu (1 - \mu) > 0$ and that $w(z, \mu)$ is linear
with respect to $z$, we have that $w(z, \mu) > 0$ for any $\mu \in
(0, 1)$ and $z \in (\mu, 3 - 2 \mu)$.  Consequently, $\f{\pa
\mscr{M}(z,\mu)}{\pa z}  < 0$ for any $\mu \in (0, 1)$ and $z \in
(\mu, 3 - 2 \mu)$.  In view of the negative sign of the partial
derivative and the fact that $\lim_{z \to 3 - 2 \mu} \mscr{M}(z,\mu)
= - \iy, \; \mscr{M}(\mu,\mu) = 0$, we have that, for any fixed $\mu
\in (0, 1)$, $\mscr{M}(z, \mu)$ is monotonically decreasing from $0$
to $- \iy$ as $z$ increases from $ \mu$ to $3 - 2 \mu$. This
completes the proof of the lemma.

\epf

The following lemma is a slight variation of Theorem 2 at page 1271
of \cite{Massart:90}, which was obtained by Massart as a byproduct
in determining the tight constant in the Dvoretzky-Kiefer-Wolfowitz
inequality.

\beL \la{lem11} Let $\ovl{X}_n = \f{\sum_{i=1}^n X_i}{n}$ where
$X_1, \cd, X_n$ are i.i.d. random variables such that $0 \leq X_i
\leq 1$ and $\bb{E}[X_i] = \mu \in (0,1)$ for $i = 1, \cd, n$. Then,
 {\small $\Pr \li \{ \ovl{X}_n \geq z
\ri \}  <  \exp \li ( n \mscr{M} (z, \mu) \ri )$} for any $z \in
(\mu, 1)$. \eeL

We can extend Lemma \ref{lem11} as follows.

 \beL \la{lem12} Let $\ovl{X}_n =
\f{\sum_{i=1}^n X_i}{n}$ where $X_1, \cd, X_n$ are i.i.d. random
variables such that $0 \leq X_i \leq 1$ and $\bb{E}[X_i] = \mu \in
(0,1)$ for $i = 1, \cd, n$. Then,
 {\small $\Pr \li \{ \ovl{X}_n \geq z
\ri \}  <  \exp \li ( n \mscr{M} (z, \mu) \ri )$} for any $z \in
(\mu, 3 - 2 \mu)$. Similarly, {\small $\Pr \li \{ \ovl{X}_n \leq z
\ri \} < \exp \li ( n \mscr{M} (z, \mu)  \ri )$} for any $z \in (- 2
\mu, \mu)$. \eeL

\bpf To show {\small $\Pr \li \{ \ovl{X}_n \geq z \ri \}  <  \exp
\li ( n \mscr{M} (z, \mu) \ri )$} for any $z \in (\mu, 3 - 2 \mu)$,
we shall consider three cases: (i) $z \in (\mu, 1)$; (ii) $z = 1$;
(iii) $z \in (1, 3 - 2 \mu)$.

In Case (i), the statement has been established as Lemma
\ref{lem11}.

In Case (ii),   we have $\Pr \li \{ \ovl{X}_n \geq z \ri \} = \Pr
\li \{ \ovl{X}_n = 1 \ri \} = \prod_{i = 1}^n \Pr \{ X_i = 1\} \leq
\prod_{i = 1}^n \bb{E} [ X_i ] = \mu^n$.  We claim that $\ln (\mu) <
\mscr{M} (1, \mu)$.  To prove this claim, it suffices to show $\ln
(\mu) < \f{ 9 (\mu - 1) } {4 \li ( 2 \mu + 1 \ri ) } $ for any $\mu
\in (0, 1)$, since $\mscr{M} (1, \mu) = \f{ 9 (\mu - 1) } {4 \li ( 2
\mu + 1 \ri ) } $.  For simplicity of notation, define {\small $g
(\mu) = \ln (\mu) - \f{ 9 (\mu - 1) } {4 \li ( 2 \mu + 1 \ri ) }$}.
Then, the first derivative of $g (\mu)$ with respect to $\mu$ is
$g^\prime (\mu) = \f{ 5 \mu^2 + 4 - 11 \mu (1 - \mu)  } { 4 \mu \li
( 2 \mu + 1 \ri )^2 } \geq \f{ 5 \mu^2 + 4 - 11 \times \f{1}{4} } {
4 \mu \li ( 2 \mu + 1 \ri )^2 } > 0$ for any $\mu \in (0, 1)$. This
implies that $g(\mu)$ is monotonically increasing with respect to
$\mu \in (0, 1)$. By virtue of such monotonicity and the fact that
$g(1) = 0$, we can conclude that $g (\mu) < 0$ for any $\mu \in (0,
1)$.  This establishes our claim that $\ln (\mu) <  \mscr{M} (1,
\mu)$.  It follows that {\small $\Pr \li \{ \ovl{X}_n \geq z \ri \}
<  \exp \li ( n \mscr{M} (z, \mu) \ri )$} holds for $z = 1$.

In Case (iii), since $0 \leq \ovl{X}_n \leq 1$,  we have $\Pr \li \{
\ovl{X}_n \geq z \ri \} = 0 < \exp \li ( n \mscr{M} (z, \mu) \ri )$
for $z \in (1, 3 - 2 \mu)$.

\bsk

To show {\small $\Pr \li \{ \ovl{X}_n \leq z \ri \}  <  \exp \li ( n
\mscr{M} (z, \mu) \ri )$} for any $z \in (- 2 \mu, \mu)$, we shall
consider three cases as follows.

In the case of $z \in (0, \mu)$,  we define $y = 1 - z$ and
$\ovl{Y}_n = \f{\sum_{i=1}^n Y_i}{n}$ with $Y_i = 1 - X_i$ for $i =
1, \cd, n$.  Then, $\Pr \li \{ \ovl{X}_n \leq z \ri \} = \Pr \li \{
\ovl{Y}_n \geq y \ri \}$. Applying Lemma \ref{lem11} to i.i.d.
random variables $Y_1, \cd, Y_n$, we have $\Pr \li \{ \ovl{Y}_n \geq
y \ri \} < \exp \li ( n \mscr{M} (y, 1 - \mu) \ri ) = \exp \li ( n
\mscr{M} (z, \mu) \ri )$ for $1 - \mu < y < 1$, i.e., $0 < z < \mu$.
This shows that {\small $\Pr \li \{ \ovl{X}_n \leq z \ri \}  <  \exp
\li ( n \mscr{M} (z, \mu) \ri )$} holds for $z \in (0, \mu)$.

In the case of $z = 0$,   we have $\Pr \li \{ \ovl{X}_n \leq z \ri
\} = \Pr \li \{ \ovl{X}_n = 0 \ri \} = \prod_{i = 1}^n (1 - \Pr \{
X_i \neq 0 \}) \leq \prod_{i = 1}^n (1 - \bb{E} [ X_i ] ) = (1 -
\mu)^n$. We claim that $\ln (1 - \mu) < \mscr{M} (0, \mu)$.  To
prove this claim, it suffices to show $\ln (1 - \mu) < \f{ 9 \mu} {4
\li ( 2 \mu - 3 \ri ) } $ for any $\mu \in (0, 1)$, since $\mscr{M}
(0, \mu) =  \f{ 9 \mu } {4 \li ( 2 \mu - 3 \ri ) } $.  For
simplicity of notation, define {\small $h (\mu) = \ln (1 - \mu) -
\f{ 9 \mu } {4 \li ( 2 \mu - 3 \ri ) }$}. Then, the first derivative
of $h (\mu)$ with respect to $\mu$ is $h^\prime (\mu) = \f{ - 16
\mu^2 + 21 \mu - 9  } { 4 (1 - \mu) \li ( 2 \mu - 3 \ri )^2 } \leq
\f{ 16 \times ( \f{21}{32})^2 - 9  } { 4 (1 - \mu) \li ( 2 \mu - 3
\ri )^2 } < 0$ for any $\mu \in (0, 1)$. This implies that $h(\mu)$
is monotonically decreasing with respect to $\mu \in (0, 1)$. By
virtue of such monotonicity and the fact that $h(0) = 0$, we can
conclude that $h (\mu) < 0$ for any $\mu \in (0, 1)$. This
establishes our claim that $\ln (1 - \mu) < \mscr{M} (0, \mu)$.  It
follows that {\small $\Pr \li \{ \ovl{X}_n \leq z \ri \} < \exp \li
( n \mscr{M} (z, \mu) \ri )$} holds for $z = 0$.

In the case of $z \in (- 2 \mu, 0)$, since $0 \leq \ovl{X}_n \leq
1$,  we have $\Pr \li \{ \ovl{X}_n \leq z \ri \} = 0 < \exp \li ( n
\mscr{M} (z, \mu) \ri )$ for $z \in (- 2 \mu, 0)$.

This completes the proof of the lemma.

\epf

Now we are in a position to prove Theorem 3.  By Lemma \ref{lemcap},
we have that, for any $\eta \in (0, 1)$, there exist two real
numbers $z_1 \in (- 2 p, p)$ and $z_2 \in (p, 3 - 2 p)$ such that
$\exp \li ( n \mscr{M} (z_1, p) \ri ) = \exp \li ( n \mscr{M} (z_2,
p) \ri ) = \f{\eta}{2}$.  Observing that $\exp \li ( n \mscr{M} (z,
p) \ri ) = \f{\eta}{2}$ can be transformed into a quadratic equation
with respect to $z$, we can obtain explicit expressions for $z_1$
and $z_2$ as {\small \[ z_1  =  p + \f{ 1 - 2 p - \sq{ 1 + \f{18 n p
(1-p)}{ \ln \f{2}{\eta}} } } {\f{2}{3} + \f{3n}{ \ln \f{2}{\eta}}},
\qqu z_2 =
 p + \f{ 1 - 2 p + \sq{ 1 +
\f{18 n p (1-p)}{ \ln \f{2}{\eta}} } } {\f{2}{3} + \f{3n}{\ln
\f{2}{\eta}}}.
\]}
Hence, by Lemma \ref{lem12}, we have
\[
\Pr \{ K \leq n z_1 \} < \exp \li ( n \mscr{M} (z_1, p) \ri ) =
\f{\eta}{2}, \qqu  \Pr \{ K \geq n z_2 \} < \exp \li ( n \mscr{M}
(z_2, p) \ri ) = \f{\eta}{2}
\]
and
\[
\mrm{T}^{-}  = \max \li \{ a, \; \li \lc n z_1  \ri \rc \ri \}, \qqu
\mrm{T}^{+} = \min \li \{ b, \; \li \lf n z_2  \ri \rf \ri \}.
\]
It follows that
\[
\Pr \{ K > \li \lf n z_2  \ri \rf \} = \Pr \{ K \geq \li \lf n z_2
\ri \rf + 1 \} \leq \Pr \{ K >  n z_2  \} \leq \Pr \{ K \geq  n z_2
\} < \f{\eta}{2},
\]
\[
\Pr \{ K < \li \lc n z_1  \ri \rc \} = \Pr \{ K \leq \li \lc n z_1
\ri \rc - 1 \} \leq \Pr \{ K <  n z_1  \} \leq \Pr \{ K \leq  n z_1
\} < \f{\eta}{2}.
\]
Since
\[
\Pr \{ a \leq K \leq b \} \leq \Pr \{ \mrm{T}^- \leq K \leq
\mrm{T}^+ \} + \Pr \{ \mrm{T}^+ <  K  \leq b \} + \Pr \{ a \leq K <
\mrm{T}^- \}
\]
and
\[
\Pr \{ \mrm{T}^+ <  K  \leq b \} \leq \Pr \{ K >  \li \lf n z_2  \ri
\rf \}, \qqu \Pr \{ a \leq K < \mrm{T}^- \} \leq \Pr \{ K < \li \lc
n z_1  \ri \rc \},
\]
we have
\[
\Pr \{ a \leq K \leq b \} < \Pr \{ \mrm{T}^- \leq K \leq \mrm{T}^+
\} + \f{\eta}{2} + \f{\eta}{2} = \Pr \{ \mrm{T}^- \leq K \leq
\mrm{T}^+ \} + \eta.
\]
On the other hand, $\Pr \{ \mrm{T}^- \leq K \leq \mrm{T}^+ \} \leq
\Pr \{ a \leq K \leq b \}$ is trivially true. This completes the
proof of Theorem 3.

\end{document}